\documentclass[10pt]{amsart}
\usepackage{amsfonts}
\usepackage{ifthen}
\usepackage{amsthm}
\usepackage{amsmath}
\usepackage{graphicx}
\usepackage{amscd,amssymb,amsthm}
\usepackage{graphicx}

\newcounter{minutes}
\divide\time by 60
\newcounter{hours}
\multiply\time by 60 \addtocounter{minutes}{-\time}

\title[von Mises--Fisher distribution]{Remarks on a parameter estimation for von Mises--Fisher distributions$^{\bigstar}$}

\author[\'Arp\'ad Baricz]{\'Arp\'ad Baricz}
\address{Department of Economics, Babe\c{s}-Bolyai University, 400591 Cluj-Napoca, Romania} \email{bariczocsi@yahoo.com}


\thanks{$^{\bigstar}$Research supported by Romanian National Research Council CNCS-UEFISCSU, project number PN-II-RU-TE\underline{ }190/2013.}

\begin{document}

\def\thefootnote{}
\footnotetext{ \texttt{File:~\jobname .tex,
          printed: \number\year-0\number\month-\number\day,
          \thehours.\ifnum\theminutes<10{0}\fi\theminutes}
} \makeatletter\def\thefootnote{\@arabic\c@footnote}\makeatother

\maketitle


\begin{abstract}
We point out an error in the proof of the main result of the paper of Tanabe et al. (2007) concerning a parameter estimation for von Mises--Fisher distributions, we correct the proof of the main result and we present a short alternative proof.
\end{abstract}

Recently, Tanabe et al. (2007) proposed an iterative algorithm using fixed points to obtain the maximum likelihood estimate for one of the parameters of the $p$-variate von Mises--Fisher distribution on the $p$-dimensional unit hypersphere. In their study Tanabe et al. (2007) arrived at the equation
\begin{equation}\label{eq1}\frac{1}{r_{\frac{p}{2}-1}(\hat{\kappa})}=\overline{R},\end{equation}
where $r_{\frac{p}{2}-1}(\hat{\kappa})=I_{\frac{p}{2}-1}(\hat{\kappa})/I_{\frac{p}{2}}(\hat{\kappa}),$ $I_{\nu}$ is the modified Bessel function of the first kind of order $\nu$, and $\overline{R}=||x_1+x_2+\dots+x_n||/n$ is the mean length of the data vector $(x_1,x_2,\dots,x_n).$ In order to solve \eqref{eq1} Tanabe et al. (2007) used a fixed point iteration method and for this first derived the following main result: if $\nu\geq 1,$ then the function $\Phi_{2\nu}:(0,\infty)\to\mathbb{R},$ defined by $\Phi_{2\nu}(x)=\overline{R}xr_{\nu-1}(x),$ has a unique fixed point.

The purpose of this note is threefold: to point out an error in the proof of above main result of Tanabe et al. (2007), to provide a correct proof and to show that the above mentioned result is almost immediate by using some known results on the ratio $1/r_{\nu-1}.$

In what follows we list our comments concerning the above mentioned result:

\begin{enumerate}
\item[\bf 1.] Let $$S_{\nu}(x)=I_{\nu}^2(x)-I_{\nu-1}(x)I_{\nu+1}(x)$$ be the Tur\'anian of the modified Bessel function of the first kind. In the proof of the above mentioned result Tanabe et al. (2007) considered the Tur\'an type inequalities
\begin{equation}\label{eq2}0<S_{\nu}(x)<\frac{1}{\nu+x}\cdot I_{\nu}^2(x)\end{equation}
and attributed these inequalities to Nasell (1974), and Thiruvenkatachar and Nanjundiah (1951). We would like to point out that there is no Tur\'an type inequality proved by Nasell (1974) and the right-hand side of \eqref{eq2} cannot be found in the paper of Thiruvenkatachar and Nanjundiah (1951). In their paper there is a Tur\'an type inequality of this kind, but with $\nu+1$ instead of $\nu+x.$ Incidentally the right-hand side of \eqref{eq2} can be found in the paper of Joshi and Bissu (1991), but it was pointed out very recently by Baricz (2012) that this inequality is not valid. This means that the right-hand side of the inequality
$$-\frac{1}{x}r_{\nu}(x)<r_{\nu}'(x)<\frac{1}{\nu+x+1}-\frac{1}{x}r_{\nu}(x)$$
is not valid, i.e. the Lemma in Tanabe et al. (2007) is not true and hence the proof of the main Theorem is not correct.

\item[\bf 2.] The proof of the main Theorem of Tanabe et al. (2007) was based on the Banach fixed point theorem, and therefore they wanted to prove that the function $\Phi_{2\nu}$ is a contraction mapping. For this it was enough to prove that $0<\Phi_{2\nu}'(x)<1$ for each $x>0$ and $\nu\geq1.$ However, since the right-hand side of \eqref{eq2} is not true, the proof of the inequality $\Phi_{2\nu}'(x)<1$ presented in Tanabe et al. (2007) is not correct. All the same, this result is true. For this observe that
    \begin{equation}\label{eq3}\Phi_{2\nu}'(x)=\overline{R}\left[r_{\nu-1}(x)+xr_{\nu-1}'(x)\right]=
    \overline{R}\left[\frac{xS_{\nu}(x)}{I_{\nu}^2(x)}\right].\end{equation}
    Since $\overline{R}\in(0,1)$ to prove $\Phi_{2\nu}'(x)<1$ it would enough to show the Tur\'an type inequality
    \begin{equation}\label{eq4}S_{\nu}(x)<\frac{1}{x}\cdot I_{\nu}^2(x).\end{equation}
    But, this inequality is valid for $\nu\geq \frac{1}{2}$ and $x>0$, as it was shown recently by Baricz (2012), and is equivalent to the fact that the function $x\mapsto xI_{\nu}'(x)/I_{\nu}(x)-x$ is strictly decreasing on $(0,\infty)$ for $\nu\geq \frac{1}{2},$ which was proved by Gronwall (1932). It is important to note here that Hamsici and Martinez (2007) used also the right-hand side of \eqref{eq2} in modeling the data of two spherical-homoscedastic von Mises-Fisher distribution. The correction of their result based on the right-hand side of \eqref{eq2} was made also by using \eqref{eq4}, see Baricz (2012) for more details. We also note that for $\nu\geq \frac{1}{2}$ and $x>0$ the Tur\'an type inequality \eqref{eq4} can be improved as (see Baricz (2012))
    $$S_{\nu}(x)<\frac{1}{\sqrt{x^2+\nu^2-\frac{1}{4}}}\cdot I_{\nu}^2(x),$$ but the contraction constant $1$ in the inequality $\Phi_{2\nu}'(x)<1$ cannot be improved. For this consider the Tur\'an type inequality (see Segura (2011))
    $$\frac{1}{\nu+\frac{1}{2}+\sqrt{x^2+\left(\nu+\frac{1}{2}\right)^2}}\cdot I_{\nu}^2(x)<S_{\nu}(x),$$
    which is valid for all $\nu\geq 0$ and $x>0.$ Combining this inequality with \eqref{eq3} and \eqref{eq4} it is easy to see that for $\nu\geq\frac{1}{2}$ the expression $xS_{\nu}(x)/I_{\nu}^2(x)$ tends to $1$ as $x$ tends to infinity. This shows that indeed the constant $1$ in the inequality $\Phi_{2\nu}'(x)<1$ cannot be improved. Now, recall that $\Phi_{2\nu}'(x)>0$ for all $\nu>0$ and $x>0,$ and combining this with the above inequality, we obtain that if $\nu\geq \frac{1}{2},$ then the function $\Phi_{2\nu}$ is a contraction mapping and therefore has a unique fixed point. This corrects the proof of the main Theorem of Tanabe et al. (2007) and improves the range of validity of the parameter $\nu=\frac{p}{2}.$

\item[\bf 3.] Finally, we present a short proof for the following result, which improves the range of validity of the main result of Tanabe et al. (2007): if $\nu>0,$ then the function $\Phi_{2\nu}$ has a unique fixed point. Observe that this statement is equivalent to the fact that the equation $1/r_{\nu-1}(x)=\overline{R}$ has a unique solution on $(0,\infty)$ for $\nu>0.$ It is known that (see Yuan and Kalbfleisch (2000)) the function $x\mapsto 1/r_{\nu-1}(x)=I_{\nu}(x)/I_{\nu-1}(x)$ is increasing on $(0,\infty)$ for $\nu\geq\frac{1}{2},$ while for $\nu\in\left(0,\frac{1}{2}\right)$ is increasing first to reach a maximum and then decreasing. Moreover, $1/r_{\nu-1}(x)\to 1$ as $x\to \infty$ for each $\nu>0,$ and the graph of $x\mapsto 1/r_{\nu-1}(x)$ approaches the asymptote from above when $\nu\in\left(0,\frac{1}{2}\right),$ and from below when $\nu\geq \frac{1}{2}.$ Now, since $\overline{R}<1,$ the above results imply that the graph of $x\mapsto 1/r_{\nu-1}(x)$ intersects only once the horizontal line $y=\overline{R}$ for each $\nu>0,$ and
    indeed the equation $1/r_{\nu-1}(x)=\overline{R}$ has a unique solution on $(0,\infty)$ for $\nu>0.$

\end{enumerate}

\end{document}